\documentclass[12pt]{article}
\usepackage{amscd,amsfonts,amssymb,amsmath,amsthm,latexsym,mathrsfs,graphicx}
\usepackage[margin=1in]{geometry}
\newcommand{\bbr}[1]{{[\![#1]\!]}}

\newcommand{\bfQ}{\mathbf{Q}}
\newcommand{\bfQp}{{\mathbf{Q}_p}}

\newcommand{\bfR}{\mathbf{R}}

\newcommand{\bfZ}{\mathbf{Z}}

\newcommand{\calC}{\mathcal{C}}
\newcommand{\calD}{\mathcal{D}}

\newcommand{\cn}{\colon}

\newcommand{\ep}{\epsilon}

\DeclareMathOperator{\Ext}{Ext}

\DeclareMathOperator{\Gr}{Gr}

\DeclareMathOperator{\Hom}{Hom}

\DeclareMathOperator{\img}{img}

\newcommand{\inv}{^{-1}}

\DeclareMathOperator{\ord}{ord}

\newcommand{\ov}[1]{{\overline{#1}}}

\newcommand{\rmP}{\mathrm{P}}

\newcommand{\rmw}{\mathrm{w}}

\newcommand{\vphi}{\varphi}

\newcommand{\wt}[1]{{\widetilde{#1}}}

\DeclareMathOperator{\rk}{rk}

\newcommand{\rmHN}{\text{HN}}
\newcommand{\wtker}{\mathop{\wt{\ker}}}

\newcommand{\wtcap}{\mathop{\wt{\cap}}}
\begin{document}

\begin{center}
{\bf Harder--Narasimhan theory}
\end{center}

These notes contain not one original idea.  An abstract formulation of
Harder--Narasimhan theory is stated without proof in [F], and I found
it helpful to write it all out.  For applications of this system, see
[F] and [K], as well as their various references; we give a rather
generic example in \S8.

\vskip 6pt

{\it (Note: I wrote this document around 2010, and was recently
  encourged to place it on the arXiv.  I've taken the opportunity to
  correct some typos, sins of exposition and formatting, and actual
  errors.  A contemporary literature search would reveal that this
  theory has since been generalized. --- Jonathan Pottharst)}

\vskip 12pt

\noindent {\bf \S1.  Hypotheses and statements.}

\vskip 6pt

The inputs and axioms given here will be in force throughout this
note.

\vskip 6pt

{\bf (1.1)} When working in exact categories, I like to say ``strict
mono/epimorphism'' as opposed to ``admissible mono/epimorphism''.  A
{\it strict subobject} is a subobject whose inclusion map is a strict
monomorphism.  I also call the distinguished (short) exact sequences
``(short) strict exact sequences''.

\vskip 6pt

\noindent \mbox{\bf Inputs (1.2).}  Our input consists of the data:
\begin{itemize}
\item a totally ordered abelian group $V$ (written additively),
\item an exact cateogry $\calC$ (our target),
\item an abelian category $\calD$,
\item a functor $F \cn \calC \to \calD$,
\item a rule $\rk \cn |\calD| \to \bfZ_{\geq 0}$ (where $|\cdot|$
  denotes isomorphism classes), and
\item a rule $\deg \cn |\calC| \to V$.
\end{itemize}

\noindent \mbox{\bf Axioms (1.3).} We assume:
\begin{itemize}
\item $F$ is exact and faithful, and for each object $X$ of $\calC$,
  $F$ induces a bijection
  \[
  \{\text{strict subobjects of }X\}
  \stackrel{\sim}{\longrightarrow}
  \{\text{subobjects of }F(X)\};
  \]
\item $\rk$ is additive over short exact sequences, and, for each
  object $X$ of $\calD$, $\rk X = 0$ if and only if $X = 0$; and
\item $\deg$ is additive over short strict exact sequences, and, for
  each morphism $f \cn X \to Y$ of $\calC$ with $F(f)$ an isomorphism,
  one has $\deg(X) \leq \deg(Y)$, with equality if and only if $f$
  itself is an isomorphism.
\end{itemize}

{\bf (1.4)} We call $V \otimes_\bfZ \bfQ$ the set of {\it slopes}.
(Typically, $V \subseteq \bfR$.)  We say that the strict subobject
$X'$ of $X$ and the subobject $F(X')$ of $F(X)$ {\it correspond} to
one another.  We call $\rk$ the {\it rank}, $\deg$ the {\it degree},
and, for nonzero objects $X$ of $\calC$, $\mu(X) := \deg(X)/\rk(X) \in
V \otimes_\bfZ \bfQ$ its {\it slope}.  (Throughout, we abusively write
$\rk(X)$ for $\rk(F(X))$ for objects $X$ of $\calC$.)  Note that the
second part of the axiom for $\deg$ can be equivalently phrased with
$\mu$ in place of $\deg$.

\vskip 6pt

{\bf (1.5)} One calls a nonzero object $X$ of $\calC$ {\it semistable}
(resp.\ {\it stable}) if for all nonzero proper subobjects $X'$ of $X$
one has $\mu(X') \leq \mu(X)$ (resp.\ $\mu(X') < \mu(X)$).

A {\it HN filtration} of $X$ is an increasing, finite, separating and
exhaustive filtration by strict subobjects
\[
0 = X_0 \subsetneq X_1 \subsetneq \cdots \subsetneq X_N = X
\]
with each $\Gr_i = X_i/X_{i-1}$ semistable, and each $\mu(\Gr_{i+1}) <
\mu(\Gr_i)$.  The $\mu(\Gr_i)$ are called the {\it HN slopes} of $X$,
and the $\rk(\Gr_i)$ are called their respective ranks.  (Note that,
in applications, simple semistable objects might have rank greater
than $1$, so that the ranks of the graded pieces can be larger than
the traditional ``multiplicities'' of the slopes.)

\vskip 6pt

{\bf (1.6)} These are the central results.

\vskip 6pt

\noindent \mbox{\bf Theorem.} {\it For any slope $\mu$, the full
  subcategory $\calC^{(\mu)}$ of $\calC$ consisting of semistable
  objects of slope $\mu$ (together with the zero object) is an
  abelian, exact subcategory of $\calC$ that is closed under
  extensions in $\calC$.}

\vskip 6pt

\noindent \mbox{\bf Theorem.} {\it Each nonzero object of $\calC$
  possesses a unique HN filtration.  When the HN filtrations are
  reindexed ``by slopes'', forming the HN filtration is functorial.}

\vskip 6pt

The first theorem is (4.12) below.  The second theorem is the
combination of (6.3), (6.2.b), and (6.5.c).

\vskip 12pt

\noindent {\bf \S2.  Properties of $F$.}

\vskip 6pt

It is impressive how much mileage one gets from only the knowledge of
$F$.

\vskip 6pt

\noindent \mbox{\bf Proposition (2.1).} {\it For any object $X$ of
  $\calC$, the correspondence induced by $F$,
  \[
  \{\text{strict subobjects of }X\}
  \stackrel{\sim}{\to}
  \{\text{subobjects of }F(X)\},
  \]
  is an isomorphism of posets.}

\vskip 6pt

\noindent {\it Proof.}  The above correspondence is a bijection by
axiom.  Let $X',X''$ be two strict subobjects of $X$.  It is clear
that if $X' \subseteq X''$ then $F(X') \subseteq F(X'')$; we show the
converse.  Consider the composite $X' \subseteq X \to X/X''$.
Applying $F$, we get $F(X') \to F(X/X'') = F(X)/F(X'')$.  If $F(X')
\subseteq F(X'')$ then this is the zero map.  Since $F$ is faithful,
$X' \to X/X''$ is the zero map, and therefore $X' \subseteq X''$.
$\square$

\vskip 6pt

{\bf (2.2)} A morphism $f \cn X \to Y$ in $\calC$ is said to be an
{\it $F$-isomorphism} (resp.\ {\it $F$-mono\-morphism}, {\it
$F$-epimorphism}) if $F(f)$ is an isomorphism (resp.\ monomorphism,
epimorphism).

\vskip 6pt

\noindent \mbox{\bf Proposition (2.3).} {\it If a morphism $f \cn X
  \to Y$ in $\calC$ is an $F$-monomorphism, then it is a monomorphism.
  In particular, if $f$ is an $F$-isomorphism, then $f$ is a
  monomorphism.}

\vskip 6pt

\noindent {\it Proof.} Suppose that $f$ is an $F$-monomorphism, and
that $g \cn Z \to X$ is a map in $\calC$ with $f \circ g = 0$.  Then
$F(f) \circ F(g) = 0$, and because $f$ is an $F$-monomorphism we have
$F(g) = 0$.  But $F$ is faithful, so $g = 0$.  $\square$

\vskip 6pt

{\bf (2.4)} Let $X$ be an object of $\calC$, and $X'$ a subobject.  We
let $\wt{X}'$ denote the strict subobject of $X$ corresponding to the
subobject $F(X')$ of $F(X)$, i.e.\ the unique strict subobject of $X$
satisfying $F(\wt{X}') = F(X')$.  We call $\wt{X}'$ the {\it                     
saturation} of $X'$ in $X$.

\vskip 6pt

\noindent \mbox{\bf Proposition (2.5).} {\it Let $X$ be an object of
  $\calC$, and $X',X''$ subobjects.  Then:
  \begin{itemize}
  \item[(a)] $X' \subseteq \wt{X}'$.
  \item[(b)] $\wt{X}' \subseteq \wt{X}''$ if and only if $F(X')
    \subseteq F(X'')$.  In particular, $F(X') = F(X'')$ if and only if
    $\wt{X}' = \wt{X}''$, and $\wt{\wt{X}'} = \wt{X}'$.
  \item[(c)] Considering $X'$ and $X''$ equivalent if $F(X') =
    F(X'')$, each equivalence class contains a unique strict object,
    equal to the saturation of any object in the class, and this
    object is the final member with respect to inclusion.
  \item[(d)] $\wt{X}'$ is characterized among strict subobjects of $X$
    containing $X'$ as being the initial one, or as the unique one
    such that $F$ takes the inclusion $X' \subseteq \wt{X}'$ to an
    isomorphism.
  \end{itemize}}

\vskip 6pt

\noindent {\it Proof.} (a): If it were true that $X' \nsubseteq
\wt{X}'$, then the natural map $X' \to X/\wt{X}'$ would be nonzero.
Since $F$ is faithful, the induced map $F(X') \to F(X/\wt{X}') =
F(X)/F(\wt{X}') = F(X)/F(X')$ would also be nonzero, but this is
impossible.

(b,c,d): These are obvious, given (a).  $\square$

\vskip 6pt

{\bf (2.6)} The process of modifying the subobject $X'$ to $\wt{X}'$,
i.e.\ one with better exactness properties, is similar to passing to a
derived or homotopy category.  From this optic, it would be tempting
to localize $\calC$ with respect to $F$, i.e.\ formally invert all
$F$-isomorphisms.  The problem is that $F$-isomorphisms do not
necessarily respect $\deg$.  In fact, by axiom, they only respect
$\deg$ when they are actual isomorphisms!

Nonetheless, one still can perform most abelian category operations in
$\calC$, as we now show.

\vskip 6pt

\noindent \mbox{\bf Proposition (2.7).}  {\it Every morphism $f \cn X
  \to Y$ has a kernel and image in $\calC$.  The kernel is a strict
  subobject of $X$, and $X/\ker f$ is canonically identified with the
  image.}

\vskip 6pt

\noindent {\it Proof.} We write, temporarily, $\wtker f$ for the
strict subobject of $X$ corresponding to $\ker F(f)$ via (2.1).  We
claim it is the kernel of $f$.  First, consider the composite $\wtker
f \to X \to Y$.  Applying $F$ gives the composite $\ker F(f) \to F(X)
\to F(Y)$, which is zero, but $F$ is faithful, so our original
composite must be zero.  On the other hand, consider any subobject $X'
\subseteq X$ such that the composite $X' \to X \to Y$ is zero.
Applying $F$, we see that $F(\wt{X}') = F(X') \to F(X) \to F(Y)$ is
zero, and therefore $F(\wt{X}') \subseteq \ker F(f)$ inside $F(X)$.
By (2.5.c), we find that $\wt{X}' \subseteq \wtker f$, and therefore
$X' \subseteq \wt{X}' \subseteq \wtker(f)$.  This shows that the
strict subobject $\wtker f$ is the kernel, so we henceforth denote it
by $\ker f$.

Now $f$ induces a map $\ov{f} \cn X/\ker f \to Y$, and applying $F$
yields the monomorphism $F(X)/\ker F(f) \to F(Y)$ in $\calD$.  By
(2.2), $\ov{f}$ is a monomorphism.  Therefore, $X/\ker f$ is
canonically identified to a subobject of $Y$ that we call $\img f$.
We claim it is an image for this morphism.  We must show that any
subobject $Y' \subseteq Y$ through which $f$ factors must contain
$\img f$.  Because we have the factorization of $f$ into $X
\stackrel{f_0}{\to} Y' \subseteq Y$, with the second map a
monomorphism, we have $\ker f = \ker f_0$.  Therefore, composite
$X/\ker f \stackrel{\ov{f}_0}{\to} Y' \subseteq Y$ is simply $\ov{f}$,
i.e.\ $\ov{f}$ and hence $f$ factors through $Y'$.  But $\ov{f}$ is
identified with the inclusion $\img f \subseteq Y$, so we are done.
$\square$

\vskip 6pt

\noindent \mbox{\bf Proposition (2.8).} {\it Let $X$ be an object of
  $\calC$, and let $X',X''$ be two subobjects.  There exists a least
  upper bound $X'+X''$ and a greatest lower bound $X' \cap X''$ for
  $X',X''$ in the poset of all subobjects of $X$, and the canonical
  sequence
  \[
  0 \to X' \cap X'' \xrightarrow{(1,-1)} X' \oplus X'' \to X' + X''
  \to 0
  \]
  is short strict exact.  If $X',X''$ are strict in $X$, then so is
  $X' \cap X''$.}

\vskip 6pt

\noindent {\it Proof.} Using (2.7), define $X'+X''$ to be the image of
the map $X' \oplus X'' \to X$.  It is immediate to verify that
$X'+X''$ is a least upper bound for $X',X''$ in $X$.

Let $K$ be the kernel of $X' \oplus X'' \to X$.  We show that the
compositions $K \to X' \oplus X'' \to X'$ and $K \to X' \oplus X'' \to
X''$ are monomorphisms.  Suppose given a map $Z \to K$ in $\calC$,
whose composite $Z \to K \to X' \oplus X'' \to X'$ is zero.  Consider
the two compositions
\begin{gather*}
K \to X' \oplus X'' \to X' \to X'+X'' \text{ and} \\
K \to X' \oplus X'' \to X'' \to X'+X''.
\end{gather*}
Their sum is the composition $K \to X' \oplus X'' \to X'+X''$, and
therefore is zero.  Because $X'' \to X'+X''$ is a monomorphism, the
composite $Z \to K \to X' \oplus X'' \to X''$ is zero.  Since $Z \to K
\to X' \oplus X''$ becomes zero followed by each of the projections to
$X',X''$, it must be zero.  But $K \to X' \oplus X''$ is a
monomorphism, and so $Z \to K$ is zero.  This shows that $K \to X'
\oplus X'' \to X'$ is a monomorphism, and the same argument shows that
$K \to X' \oplus X'' \to X''$ is a monomorphism as well.

We denote by $X' \cap X''$ the subobject of $X'+X''$ determined by the
monomorphism $K \to X' \oplus X'' \to X' \to X'+X''$.  By construction
it is a subobject of $X'$, and the above argument shows that it is
also a subobject of $X''$.  The same reasoning shows that the
composition $X' \cap X'' \to X' \oplus X'' \to X''$ is the {\it
negative} of the canonical inclusion $X' \cap X'' \subseteq X''$ as
subobjects of $X'+X''$.  Therefore, we have the short exact sequence
as claimed.

We now show that the subobject $X' \cap X''$ is a greatest lower bound
for $X',X''$ in $X$.  Given any subobject $Y$ of $X$ whose inclusion
into $X$ is factored through both $X',X''$, we must factor its
inclusion into $X$ through $X' \cap X''$.  Combining the two given
factorizations gives a map $(1,-1) \cn Y \to X' \oplus X''$, whose
projection to $X'$ is the inclusion $Y \subseteq X'$ and whose sum to
$X$ is zero.  Therefore, it uniquely factors through $X' \cap X'' \to
X' \oplus X''$.  But the composite $X' \cap X'' \to X' \oplus X'' \to
X'$ is the canonical inclusion, so we have exhibited a factorization
of $Y \subseteq X'$ through $X' \cap X'' \subseteq X'$.

Finally, suppose that $X',X''$ are strict.  Denote by $X' \wtcap X''$
the strict subobject of $X$ corresponding via (2.1) to $F(X') \cap
F(X'')$.  We show that the strict object $X' \wtcap X''$ is a greatest
lower bound for $X,X''$ in $X$.  By the isomorphism of posets, $X'
\wtcap X''$ is contained in both $X',X''$.  Moreover, if $Y$ is any
subobject of $X$ that is contained in $X',X''$, then $F(Y) \subseteq
F(X')$ and $F(Y) \subseteq F(X'')$, hence $F(Y) \subseteq F(X') \cap
F(X'')$ in $F(X)$.  By (2.5.b) this is equivalent to $\wt{Y} \subseteq
X' \wtcap X''$, and we conclude by (2.5.a).  $\square$

\vskip 6pt

{\bf (2.9)} We warn the reader that, in the above lemma, even if
$X',X''$ are strict subobjects of $X$, the subobject $X'+X''$ might
not be strict.

\vskip 6pt

\noindent \mbox{\bf Proposition (2.10).} {\it Let $X$ be an object of
  $\calC$, and $X',X''$ two subobjects.  Assume the inclusions $X'
  \cap X'' \subseteq X'$ and $X'' \subseteq X'+X''$ are strict.  (This
  holds, in particular, if $X',X''$ are strict subobjects of $X$.)
  Then the canonical map
  \[
  f \cn \frac{X'}{X' \cap X''} \to \frac{X'+X''}{X''}
  \]
  induced by $X' \subseteq X'+X''$ is an isomorphism.  Moreover, these
  objects are canonically subobjects of $X/X''$.}

\vskip 6pt

\noindent {\it Proof.} Applying $F$ to $f$ clearly yields an
isomorphism, so by (2.2) $f$ is a monomorphism.  We must show that the
subobject thus obtained is the whole thing.  Considering that the
compositions $X' \to X'/(X' \cap X'') \to (X'+X'')/X''$ and $X' \to
X'+X'' \to (X'+X'')/X''$ are equal, we see that the image of $f$
contains the image of $X'$ under the latter of these two.  The image
of $f$ contains $0$, which is the image of the composite $X'' \to
X'+X'' \to X'+X''/X''$.  Therefore, the image of $f$ contains the
image of the composite $X' \oplus X'' \to X'+X' \to (X'+X')/X''$.  But
each of the latter two maps are strict surjections.  The second claim
follows from applying (2.3) to the canonical map $(X'+X'')/X'' \to
X/X''$.  $\square$

\vskip 12pt

\noindent {\bf \S3.  Properties of $\rk$.}

\vskip 6pt

The main use of $\rk$ is to bound chains of strict subobjects in
$\calC$.

\vskip 6pt

\noindent \mbox{\bf Proposition (3.1).} {\it Every object of $\calD$
  has finite length.  For any object $\calD$, any nonempty collection
  of its subobjects closed under $+$ admits a final element.}

\vskip 6pt

\noindent {\it Proof.} By d\'evissage and induction, one sees that
length of an object $X$ of $\calD$ is bounded above by its rank.

For the second claim, the ranks in the collection are bounded above by
the rank of the parent object, and any member whose rank is maximal is
the final element.  $\square$

\vskip 6pt

\noindent \mbox{\bf Proposition (3.2).} {\it Let $X$ be an object of
  $\calC$, and $X' \subseteq X'' \subseteq X$ subobjects.  Then
  $\rk(X') = \rk(\wt{X}')$, and if $\rk(X') = \rk(X'')$ then $\wt{X}'
  = \wt{X}''$.  In particular, any subobject $X'$ of $X$ with $\rk(X')
  = \rk(X)$ has $\wt{X}' = X$.}

\vskip 6pt

\noindent {\it Proof.}  The first claim is clear, because $F(X') =
F(\wt{X}')$.  For the second claim, we may therefore replace $X'$ and
$X''$ by $\wt{X}'$ and $\wt{X}''$, and assume that $X'$ and $X''$ are
strict.  Now consider the short strict exact sequence $0 \to X' \to
X'' \to X''/X' \to 0$.  Applying $F$, we get the short exact sequence
$0 \to F(X') \to F(X'') \to F(X''/X') \to 0$.  We must have
$\rk(F(X''/X')) = 0$, which forces $0 = F(X''/X') = F(X'')/F(X')$.
Therefore $F(X') = F(X'')$, and the axiom for $F$ shows that $X' =
X''$ because $X'$ and $X''$ are strict.  $\square$

\vskip 12pt

\noindent {\bf \S4.  Properties of $\deg$ and $\mu$.}

\vskip 6pt

With notions of $\deg$ and $\mu$, an arithmetic of subobjects
emerges.

\vskip 6pt

\noindent \mbox{\bf Proposition (4.1).} {\it Let $X$ be an object of
  $\calC$, and consider the equivalence relation on subobjects of $X$
  from (2.5.c).  In any equivalence class, all members have the same
  rank, and the unique strict subobject is the unique member with
  maximal degree (resp.\ slope).}

\vskip 6pt

\noindent {\it Proof.} The claim for $\rk$ is trivial.  In a given
equivalence class, given that the strict object is final with respect
to inclusion, the claim for $\deg$ is a restatement of our axiom for
$\deg$.  The claim for $\mu$ follows.  $\square$

\vskip 6pt

\noindent \mbox{\bf Proposition (4.2).} {\it Let $X$ be an object of
  $\calC$, and $X'$ a subobject whose slope is maximal among the
  slopes of subobjects of $X$.  Then $X'$ is a strict subobject.}

\vskip 6pt

\noindent {\it Proof.}  By (2.5.d), one has $X' = \wt{X}'$ if and only
if $X'$ is strict.  Moreover, by (4.1) one has $\mu(X') \leq
\mu(\wt{X}')$, with equality if and only if $X' = \wt{X}'$.  But the
slope of $X'$ is maximal, so also $\mu(\wt{X}') \leq \mu(X')$.
$\square$

\vskip 6pt

\noindent \mbox{\bf Proposition (4.3).} {\it Let $0 = X_0 \subsetneq
  X_1 \subsetneq \cdots \subsetneq X_N = X$ be strict subobjects of
  $X$.  Write $\Gr_i = X_i/X_{i-1}$, $\rk_i = \rk(\Gr_i)$, $\deg_i =
  \deg(\Gr_i)$, and $\mu_i = \mu(\Gr_i)$.  Then:
  \begin{itemize}
  \item[(a)] One has the identity
  \[
  \mu(X) = \sum_{i=1}^N \frac{\rk_i}{\rk(X)}\,\mu_i.
  \]
  \item[(b)] Either $\min\, \{\mu_i\}_{i=1\ldots N} < \mu(X) <
    \max\,\{\mu_i\}_{i=1\ldots N}$ holds, or all $\mu_i$ are equal to
    $\mu(X)$.
  \end{itemize}}

\vskip 6pt

\noindent {\it Proof.}  (a): By induction, using the additivity of
$\deg$ over short exact sequences, we find that $\deg(X) = \sum
\deg_i$.  Then we calculate:
\[
\mu(X) = \frac{\deg(X)}{\mu(X)}
= \frac{\sum \deg_i}{\rk(X)}
= \frac{\sum \rk_i\mu_i}{\rk(X)}
= \sum \frac{\rk_i}{\rk(X)}\,\mu_i.
\]

(b): By induction, using the additivity of $\rk$ over short exact
sequences, we find that $\rk(X) = \sum \rk_i$.  Therefore, each of the
rational numbers $\rk_i/\rk(X)$ lies strictly between $0$ and $1$, and
their sum is $1$.  Thus, (a) expresses $\mu(X)$ as a weighted average
of the $\mu_i$ with weights strictly between $0$ and $1$ that sum to
$1$.  The result follows.  $\square$

\vskip 6pt

{\bf (4.4)} We will most often use (4.3.b) as follows: given a short
strict exact sequence $0 \to X' \to X \to X'' \to 0$, exactly one of
$\mu(X') < \mu(X) < \mu(X'')$, $\mu(X'') < \mu(X) < \mu(X')$, or
$\mu(X') = \mu(X) = \mu(X'')$ holds.  Therefore, to know one of the
ordering relationships among these three slopes, it suffices to know
either of the other two.

\vskip 6pt

\noindent \mbox{\bf Proposition (4.5).} {\it For a nonzero subobject
  $X$ of $\calC$, the following conditions are equivalent:
  \begin{itemize}
  \item For all nonzero proper subobjects $X'$ of $X$, one has
    $\mu(X') \leq \mu(X)$ (resp.\ $\mu(X') < \mu(X)$).
  \item For all nonzero, proper, strict subobjects $X'$ of $X$, one
    has $\mu(X') \leq \mu(X)$ (resp.\ $\mu(X') < \mu(X)$).
  \item For all nonzero, proper, strict subobjects $X'$ of $X$', one
    has $\mu(X) \leq \mu(X/X')$ (resp. $\mu(X) < \mu(X/X')$).
  \end{itemize}}

\vskip 6pt

\noindent {\it Proof.}  The first two conditions are equivalent by
(4.1).  The last two conditions are equivalent by (4.4).  $\square$

\vskip 6pt

{\bf (4.6)} If the equivalent conditions of (4.5) hold for $X$, we say
that $X$ is {\it semistable}.

\vskip 6pt

\noindent \mbox{\bf Proposition (4.7).} {\it Let $X$ be an object of
  $\calC$, and $X'$ a nonzero subobject whose slope is maximal among
  the slopes of nonzero subobjects of $X$.  Then $X'$ is semistable.

  In particular, if $X$ has rank one, then $X$ is semistable.}

\vskip 6pt

\noindent {\it Proof.}  The first claim follows directly from the
definition of semistability, since any subobject of $X'$ is a
subobject of $X$.  For the second claim, assume $X$ has rank one.  By
additivity and nonnegativity, any nonzero subobject of $X$ has rank
one, and by (3.2) it follows that the only nonzero strict subobject of
$X$ is $X$ itself.  By definition, then, $X$ is semistable.  $\square$

\vskip 6pt

\noindent \mbox{\bf Proposition (4.8).} {\it Let $0 \to X' \to X \to
  X'' \to 0$ be a short strict exact sequence in $\calC$.  If $X'$ and
  $X''$ are both semistable of the same slope $\mu$, then $X$ is also
  semistable of slope $\mu$.  In particular, $X' \oplus X''$ is
  semistable of slope $\mu$.}

\vskip 6pt

\noindent {\it Proof.} By (4.4), $\mu(X) = \mu$.  Let $Y \subseteq X$
be any strict subobject.  Consider the short strict exact sequence
\[
0 \to Y \cap X' \to Y \to Y/(Y \cap X') \to 0.
\]
By (4.4), we must show that $\mu(Y \cap X') \leq \mu$ and $\mu(Y/(Y
\cap X')) \leq \mu$.  The first follows because $X'$ is semistable of
slope $\mu$.  The second follows because $Y/(Y \cap X')$ is a
subobject of $X''$ by (2.10) and because $X''$ is
semistable of slope $\mu$.  $\square$

\vskip 6pt

\noindent \mbox{\bf Proposition (4.9).} {\it Let $X$ be an object of
  $\calC$ that is semistable of slope $\mu$.  Then any nonzero strict
  subobject or strict quotient of $X$ that has slope $\mu$ is
  semistable.}

\vskip 6pt

\noindent {\it Proof.} Let $X' \subseteq X$ be nonzero, strict, and of
slope $\mu$.  Then any subobject of $X'$ is also a subobject of $X$,
and therefore has slope at most $\mu$.  Therefore $X'$ is semistable.
Dually, let $X \twoheadrightarrow X''$ be a nonzero, strict, and of
slope $\mu$.  Then any strict quotient of $X''$ is also a strict
quotient of $X$, and therefore has slope at least $\mu$.  Therefore
$X''$ is semistable.  $\square$

\vskip 6pt

\noindent \mbox{\bf Proposition (4.10).} {\it Let $f \cn X \to Y$ be a
  morphism in $\calC$ with $X$ and $Y$ semistable of the same slope
  $\mu$.  Then $\ker f$ and $\img f$ are semistable of slope $\mu$,
  and $\img f$ is a strict subobject of $Y$.}

\vskip 6pt

\noindent {\it Proof.} Since $X$ (resp.\ $Y$) is semistable of slope
$\mu$, one has $\mu(\ker f) \leq \mu$ (resp.\ $\mu(\img f) \leq \mu$).
But now (4.4) and the short strict exact sequence
\[
0 \to \ker f \to X \to \img f \to 0
\]
imply that $\mu(\ker f) = \mu(\img f) = \mu$.  By (4.9), $\ker f$ and
$\img f$ are semistable as well.  Finally, by (4.2), $\img f \subseteq
Y$ is strict if $\mu(\img f)$ is maximal among the slopes of subobjects
of $Y$.  Since $Y$ is semistable of slope $\mu = \mu(\img f)$, this
condition holds.  $\square$

\vskip 6pt

{\bf (4.11)} Let $\mu$ be a slope.  We denote by $\calC^{(\mu)}$ the
full subcategory of $\calC$ consisting of semistable objects of slope
$\mu$, together with the zero object.

\vskip 6pt

\noindent \mbox{\bf Proposition (4.12).} {\it For every slope $\mu$,
  the category $\calC^{(\mu)}$ is an abelian, exact subcategory of
  $\calC$ that is closed under extensions in $\calC$.}

\vskip 6pt

\noindent {\it Proof.} By (4.10), morphisms in $\calC^{(\mu)}$ admit
kernels and images.  By (2.7) and (4.10), images are also coimages.
The existence of cokernels follows.  Hence it is abelian.  By
definition, a sequence is exact in $\calC^{(\mu)}$ if and only if it
is exact in $\calC$; this gives exactness.  The final claim is (4.8).
$\square$

\vskip 12pt

\noindent {\bf \S5.  Subobjects of maximal slope.}

\vskip 6pt

This preparatory section contains the technical arguments that make
HN filtrations work.

\vskip 6pt

{\bf (5.1)} Let $X$ be an object of $\calC$.  We say a subobject $X'$
of $X$ is {\it of maximal slope} if it is nonzero and $\mu(X'') \leq
\mu(X')$ for all nonzero subobjects $X''$ of $X$.  For example, $X$ is
semistable if and only if $X$ is itself a subobject of $X$ of maximal
slope.  Recall that by (4.7) such a subobject $X'$ is semistable, and
by (4.2) it is a strict subobject of $X$.

\vskip 6pt

\noindent \mbox{\bf Proposition (5.2).} {\it Let $X$ be an object of
  $\calC$.  Let $X',X''$ be subobjects of $X$ of maximal slope.  Then
  $X'+X''$ and $X' \cap X''$ are of maximal slope if nonzero.

  In particular, if $X$ admits any subobjects of maximal slope, then
  $X$ admits a unique final subobject of maximal slope.}

\vskip 6pt

\noindent {\it Proof.}  By (4.2), $X',X''$ are strict.  Since $X'$ is
semistable by (4.7), we have $\mu(X') \leq \mu(X'/(X' \cap X'')) =
\mu((X'+X'')/X'')$ by (2.10).  Since $\mu(X'')$ is maximal, we have
$\mu(X'+X'') \leq \mu(X'')$, and by (4.4) we have $\mu((X'+X'')/X'')
\leq \mu(X'+X'')$.  Finally, since $\mu(X')$ is maximal, we have
$\mu(X'+X'') \leq \mu(X')$.  Summarizing, we have the inequalities
\[
\mu(X') \leq \mu(X'/(X' \cap X'')) = \mu((X'+X'')/X'') \leq
\mu(X'+X'') \leq \mu(X'),
\]
which must all be equalities.  The claim for $X'+X''$ immediately
follows.  On the other hand, when $X' \cap X''$ is nonzero, knowing
that the first of these inequalities is an equality shows by (4.4)
that we must also have $\mu(X' \cap X'') = \mu(X')$, whence the claim
for $X' \cap X''$.  The final claim follows from (2.1) and (3.1).
$\square$

\vskip 6pt

{\bf (5.3)} Given an object $X$ of $\calC$, consider the condition on
a nonzero subobject $X'$ of $X$:
\[
(5.3.1) \qquad
\text{for all subobjects $X''$ of $X$ properly containing $X'$, one
  has $\mu(X'') < \mu(X')$.}
\]
By (4.1) such an $X'$ is strict, and also it suffices to know
(5.3.1) for strict $X''$.

\vskip 6pt

\noindent \mbox{\bf Proposition (5.4).} {\it Let $X$ be a an object of
  $\calC$, with nonzero strict subobjects $X',X''$.  Suppose $X'$
  satisfies (5.3.1), $X''$ is semistable, and $X'' \nsubseteq X'$.
  Then $\mu(X'') < \mu(X')$.}

\vskip 6pt

\noindent {\it Proof.} Let $f \cn X'' \to X \to X/X'$ be the
composition.  By hypothesis, $f$ is nonzero, and hence $F(f)$ is
nonzero.  Because $X''$ is semistable we have $\mu(X'') \leq \mu(\img
f)$.  By (4.1), we have $\mu(\img f) \leq \mu(\wt{\img f})$.  Let
$X^3$ be the unique strict subobject of $X$ (properly) containing $X'$
with $X^3/X' = \wt{\img f} \subseteq X/X'$.  By (5.3.1), $\mu(X^3) <
\mu(X')$.  By (4.4), we must also have $\mu(\wt{\img f}) < \mu(X^3)$.
Putting everything together, we have $\mu(X'') \leq \mu(\img f) \leq
\mu(\wt{\img f}) < \mu(X^3) < \mu(X')$.  $\square$

\vskip 6pt

\noindent \mbox{\bf Proposition (5.5).} {\it Let $X$ be an object of
  $\calC$.  The following conditions are equivalent for a nonzero
  subobject $X'$ of $X$:
  \begin{itemize}
  \item[(a)] $X'$ satisfies (5.3.1) and is semistable.
  \item[(b)] $X'$ satisfies (5.3.1) and is of maximal slope.
  \item[(c)] $X'$ is the final subobject of $X$ of maximal slope.
\end{itemize}}

\vskip 6pt

\noindent {\it Proof.} (a) $\implies$ (b): Let $X''$ be a nonzero
subobject of $X$.  To show that $\mu(X'') \leq \mu(X')$, we may
replace $X''$ by $\wt{X}''$ and assume that $X''$ is a strict
subobject.  Since $X''+X'$ contains $X'$, by (5.3.1) we have
$\mu(X''+X') \leq \mu(X')$.  Since $X'$ is semistable, by (2.10) we
have $\mu(X') \leq \mu(X'/(X'' \cap X')) = \mu((X''+X')/X'')$.
Combining these, we have $\mu(X''+X') \leq \mu((X''+X')/X'')$, and now
(4.4) implies that $\mu(X'') \leq \mu(X''+X')$.  To conclude, we
combine the preceding inequality with another application of (5.3.1),
$\mu(X''+X') \leq \mu(X')$.

(b) $\implies$ (a): This is (4.7).

(a,b) $\implies$ (c): We must show that any nonzero subobject $X''$ of
$X$ having maximal slope is contained in $X'$.  But this follows from
the contrapositive (5.4), because (5.3.1) holds for $X'$, (4.7) shows
$X''$ is semistable, and both $X',X''$ are of maximal slope so
$\mu(X'')=\mu(X')$.

(c) $\implies$ (a,b): We must check (5.3.1).  Let $X''$ be a subobject
of $X$ properly containing $X'$.  If it were possible that $\mu(X')
\leq \mu(X'')$, then because $\mu(X')$ is maximal we would have
$\mu(X'') = \mu(X')$.  But $X'$ is final among objects with slope
$\mu(X')$, so this is impossible.  $\square$

\vskip 6pt

{\bf (5.6)} Let $X$ be an object of $\calC$, and $X'$ a subobject
satisfying the equivalent conditions of (5.4).  By (5.4.c), such an
$X'$ is uniquely determined, if it exists.  By (4.2), $X'$ is a strict
subobject of $X$.  It is easy to see that $X$ is semistable if and
only if $X' = X$.  When $X' \subsetneq X$, by (5.4.c) $X'$ is the
maximal subobject that maximally contradicts the semistability of $X$.
We call $X'$ the {\it SCSS} for $X$, where ``SCSS'' stands for {\it
  strongly contradicting semistability}.

\vskip 6pt

\noindent \mbox{\bf Proposition (5.7).} {\it Every nonzero object $X$
  of $\calC$ admits an SCSS.}

\vskip 6pt

\noindent {\it Proof.}  We proceed by induction on the rank $d$ of
$X$.  The base case is where $X$ is semistable (which includes when
$d=1$ by (4.7)) is trivial, so we assume that $d>1$ and $X$ is not
semistable.  Because $X$ is not semistable, there exist nonzero,
proper submodules $Y$ of $X$ having $\mu(X) < \mu(Y)$.  If we were to
have $\rk(Y)=\rk(X)$, then by (3.2) we would have $\wt{Y}=X$, and by
(4.1) $\mu(Y) \leq \mu(X)$.  Thus, any such $Y$ must have rank
strictly less than $d$.

Therefore, there is a maximal integer $s$ among the ranks of nonzero
proper submodules $Y$ of $X$ satisfying $\mu(X) < \mu(Y)$, and $s$
satisfies $0<s<d$.  Fix any such submodule $Y$ of rank $s$.  Replacing
$Y$ by $\wt{Y}$ (which neither decreases $\mu(Y)$ nor changes
$\rk(Y)$), we may assume $Y \subseteq X$ is strict.  Apply the
inductive hypothesis to produce an SCSS $Y'$ of $Y$.  We verify that
$Y'$ is an SCSS for $X$, by showing that any nonzero strict subobject
$X'$ of $X$ satisfying $\mu(Y') \leq \mu(X')$ is contained in $Y'$.

We have the short strict exact sequence by (2.8):
\[
0 \to Y \cap X' \xrightarrow{(1,-1)} Y \oplus X' \to Y+X' \to 0.
\]
Therefore, using $\mu(Y \cap X') \leq \mu(Y') \leq \mu(X')$,
followed by $\mu(Y) \leq \mu(Y') \leq \mu(X')$, we have
\begin{align*}
\rk(Y+X')\mu(Y+X')
  &= \deg(Y+X') = \deg(Y \oplus X') - \deg(Y \cap X')\\
  &= \deg(Y) + \deg(X') - \deg(Y \cap X')\\
  &= \rk(Y)\mu(Y) + \rk(X')\mu(X') - \rk(Y \cap X')\mu(Y \cap X') \\
  &\geq \rk(Y)\mu(Y) + [\rk(X')-\rk(Y \cap X')]\mu(X') \\
  &\geq [\rk(Y) + \rk(X') - \rk(Y \cap X')]\mu(Y) \\
  &= \rk(Y+X')\mu(Y).
\end{align*}
But $\rk(Y+X') \geq \rk(Y) = s > 0$, so $\mu(Y+X') \geq \mu(Y) >
\mu(X)$.  The maximality of $s$ therefore implies that $\rk(Y+X') =
\rk(Y)$.  Now we have
\[
Y \subseteq Y+X' \subseteq \wt{Y+X'} = \wt{Y} = Y,
\]
with the the penultimate identity coming from (3.2) and the final
identity because $Y$ is strict.  But this forces $X' \subseteq Y$,
and, $Y'$ being an SCSS for $Y$, $\mu(Y') \leq \mu(X')$ forces $X'
\subseteq Y'$.  $\square$

\vskip 6pt

\noindent \mbox{\bf Acknowledgment (5.8).}  We thank Ruochuan Liu for
sharing the above argument with us.

\vskip 12pt

\noindent {\bf \S6.  HN filtrations.}

\vskip 6pt

We now assemble the preceding facts to show the main results on HN
filtrations: existence, uniqueness, and functoriality.

\vskip 6pt

{\bf (6.1)} Let $X$ be a nonzero object of $\calC$.  By an {\it HN
  filtration} of $X$ we mean an increasing, finite, separating and
exhaustive filtration by strict subobjects
\[
0 = X_0 \subsetneq X_1 \subsetneq \cdots \subsetneq X_N = X
\]
with each $\Gr_i = X_i/X_{i-1}$ semistable, and each $\mu(\Gr_{i+1}) <
\mu(\Gr_i)$.  Showing the existence and uniqueness of the HN
filtration is now straightforward, because the SCSS, studied in the
preceding section, is its unique first step.

\vskip 6pt

\noindent \mbox{\bf Proposition (6.2).} {\it Let $X$ be a nonzero
  object of $\calC$.  Then:
\begin{itemize}
\item[(a)] If $X$ is not semisimple, and $0 = X_0 \subsetneq X_1
  \subsetneq \cdots \subsetneq X_N = X$ is any HN filtration for $X$,
  then $X_1$ satisfies the equivalent conditions of (5.5) for an SCSS
  for $X$.
\item[(b)] If a HN filtration for $X$ exists, then it is unique.
\end{itemize}}

\vskip 6pt

\noindent {\it Proof.} (a): We induct on the length $N$ of the
filtration.  In the case $N=1$, $X$ is semistable, and we are done.
Now assume $N>1$.  By hypothesis, $X_1$ is semistable.  Given a strict
subobject $X''$ of $X$ properly containing $X_1$, we must show that
$\mu(X'') < \mu(X_1)$.  Equivalently by (4.4), we must show that
$\mu(X''/X_1) < \mu(X_1)$.  Since $\mu(\Gr_2) < \mu(\Gr_1) =
\mu(X_1)$, it suffices to show that any subobject of $X/X_1$ has slope
at most $\mu(\Gr_2)$.  But by inductive hypothesis, $\Gr_2$ satisfies
the equivalent conditions of (5.5) for $X/X_1$, and in particular
$\Gr_2$ is a subobject of $X/X_1$ of maximal slope.

(b): By (a), for each $1 \leq i < N$, $X_i/X_{i-1}$ is the SCSS of
$X/X_{i-1}$, and then $X_i$ is uniquely determined as the preimage of
$X_i/X_{i-1}$ in $X$.  $\square$

\vskip 6pt

\noindent \mbox{\bf Proposition (6.3).} {\it Each nonzero object $X$
  of $\calC$ possesses a HN filtration.}

\vskip 6pt

\noindent {\it Proof.} We perform induction on the rank of $X$.  The
base case is when $X$ is semistable, which includes the case of rank
one by (4.7), and we take $X_0 = 0$ and $X_1 = X$.  Now suppose given
$X$, assumed not semistable.  Let $X'$ be the SCSS of $X$, which
exists by (5.7).  By inductive hypothesis, $X/X'$ has a HN filtration,
$0 = Y_0 \subsetneq Y_1 \subsetneq \cdots \subsetneq Y_N = X/X'$.
Take $X_0 = 0$, and for $1 \leq i \leq N+1$ let $X_i$ be the preimage
of $Y_{i-1}$ in $X$.  By construction, all the graded pieces are
semistable.  We only need to show that $\mu(\Gr_2) < \mu(\Gr_1)$.
Considering the short strict exact sequence $0 \to \Gr_1 \to X_2 \to
\Gr_2 \to 0$, by (4.4) it suffices to show that $\mu(X_2) <
\mu(\Gr_1)$, and this follows because $\Gr_1 = X'$ satisfies (5.3.1).
$\square$

\vskip 6pt

{\bf (6.4)} Each nonzero object $X$ of $\calC$ has canonically
associated to it a collection of slopes with ranks, called the {\it HN
  slopes} and their respective {\it ranks}: if $0 = X_0 \subsetneq X_1
\subsetneq \cdots \subsetneq X_N = X$ is the HN filtration, then we
set $\Gr_i = \Gr_i(X) = X_i/X_{i-1}$, $\mu_i = \mu_i(X) = \mu(\Gr_i)$,
and $\rk_i = \rk_i(X) = \rk(\Gr_i)$.  However, one only gets good
functorial properties from using {\it slopes} to index the filtration.
For any slope $\mu$, let $X^{(\mu)}$ be the strict subobject of $X$
defined by
\[
X^{(\mu)} =
\begin{cases}
X = X_N & \text{if }\mu \leq \mu_N, \\
X_{N-1} & \text{if }\mu_N < \mu \leq \mu_{N-1}, \\
\vdots & \vdots \\
X_1 & \text{if }\mu_2 < \mu \leq \mu_1, \text{
  and}\\
0 & \text{if }\mu > \mu_1.
\end{cases}
\]
The $X^{(\mu)}$ define a decreasing, separating and exhaustive
filtration indexed by the slopes.  It has finitely many breaks,
occurring precisely at the slopes of $X$.  We let $\Gr^{(\mu)} =
\Gr^{(\mu)}(X) = X^{(\mu)}/X^{(\mu+)}$, where $X^{(\mu+)}$ is the next
smaller step in the HN filtration of $X$ if $\mu$ is a break, and
$X^{(\mu)}$ itself otherwise.  Thus, the nonzero $\Gr^{(\mu)}$ are
precisely the $\Gr_i$, indexed by their slopes: $\Gr^{(\mu_i)} =
\Gr_i$.

\vskip 6pt

\noindent \mbox{\bf Lemma (6.5).} {\it Let $X$ and $Y$ be two nonzero
  objects of $\calC$.  Then:
\begin{itemize}
\item[(a)] If $X$ and $Y$ are semistable and $\Hom_\calC(X,Y) \neq 0$,
  then $\mu(X) \leq \mu(Y)$.
\item[(b)] If all HN slopes of $X$ are strictly greater than those of
  $Y$, then $\Hom_\calC(X,Y) = 0$.
\item[(c)] Let $f \in \Hom_\calC(X,Y)$, and $\mu$ be a slope.  The
  restriction of $f$ to $X^{(\mu)}$ factors canonically through
  $Y^{(\mu)}$.
\item[(d)] For each slope $\mu$, the rule $X \mapsto \Gr^{(\mu)}(X)$
  gives a functor $\calC \to \calC^{(\mu)}$.
\end{itemize}}

\vskip 6pt

\noindent {\it Proof.} (a): Suppose $f \in \Hom_\calC(X,Y)$ is
nonzero.  As $X$ is semistable, one has $\mu(X) \leq \mu(\img f)$, and
since $Y$ is semistable, one has $\mu(\img f) \leq \mu(Y)$.

(b): Let $X$ (resp.\ $Y$) have HN filtration $0 = X_0 \subsetneq X_1
\subsetneq \cdots \subsetneq X_N = X$ (resp.\ $0 = Y_0 \subsetneq Y_1
\subsetneq \cdots \subsetneq Y_M = Y$).  We show by induction on $N$
that any map $X \to Y$ must vanish.  The base case, when $X$ is
semistable, is proved by induction on $M$.  The base case on $M$, when
$Y$ is also semistable, is a restatement of (a).  For the inductive
step on $M$, given a map $X \to Y$ we apply the base case to see that
the composite $X \to Y \to Y/Y_{M-1} = \Gr_M(Y)$ vanishes, hence the
map $X \to Y$ factors through a map $X \to Y_{M-1}$, which also must
vanish by inductive hypothesis.  For the inductive step on $N$, we
apply the base case to the composite $X_1 \to X \to Y$ to see that our
given map $X \to Y$ factors through a map $X/X_1 \to Y$, which in turn
vanishes by inductive hypothesis.

(c): Let $f$ and $\mu$ be given.  It suffices to replace $X$ by
$X^{(\mu)}$, $f$ by its restriction to $X^{(\mu)}$, and $Y$ by
$Y/Y^{(\mu)}$, whence $X = X^{(\mu)}$ and $Y^{(\mu)} = 0$, and to show
that $f = 0$.  But now we are in the situation of (b), so indeed $f =
0$.

(d): This follows trivially from (c).  $\square$

\vskip 6pt

{\bf (6.6)} Although (6.5) gives us, for any map $f \cn X \to Y$ in
$\calC$, a map $f^{(\mu)} \cn \Gr^{(\mu)}(X) \to \Gr^{(\mu)}(Y)$ in
$\calC^{(\mu)}$ for each slope $\mu$, passing from $f$ to the family
of $f^{(\mu)}$ is very destructive in practice.  For example, if $X$
has rank one, and $\iota \cn X' \subseteq X$ is a nonzero proper
subobject, then all the $\iota^{(\mu)}$ vanish.

\vskip 12pt

\noindent {\bf \S7.  HN polygons.}

\vskip 6pt

We give some complements on HN polygons, which are an analogue of
Newton polygons in this setting.

\vskip 6pt

{\bf (7.1)} Let $a,b \in \bfQ$.  By $[a,b]_\bfQ$, we mean the interval
$[a,b] \cap \bfQ$ of rational numbers.  The set of slopes $V
\otimes_\bfZ \bfQ$ is a $\bfQ$-vector space, so there is an obvious
notion of a linear function on $[a,b]_\bfQ$ with values in $V
\otimes_\bfZ \bfQ$.  Similarly, we may speak of piecewise-linear
functions on $[a,b]_\bfQ$ with values in $V \otimes_\bfZ \bfQ$, and we
understand the breakpoints to be at rational numbers.

\vskip 6pt

{\bf (7.2)} Let $X$ be a nonzero object of $\calC$, and $0 = X_0
\subsetneq X_1 \subsetneq \cdots \subsetneq X_N = X$ its HN
filtration.  We define $\rmHN(X)$ to be the unique concave-down,
piecewise-linear function defined on the rational interval
$[0,\rk(X)]_\bfQ$, satisfying $\rmHN(X)(0) = 0$, and whose $i$th slope
is $\mu_i$ with horizontal width $\rk_i$.  The {\it HN polygon} of $X$
is defined to be the graph of $\rmHN(X)$.  We call the point
$(\rk(X),\deg(X))$ the {\it HN endpoint} of $X$; it is easy to
calculate that $\rmHN(X)(\rk(X)) = \deg(X)$, thus explaining the
terminology.

Note that $0 = X_0 \subsetneq X_1 \subsetneq \cdots \subsetneq X_i =
X_i$ is the HN filtration of $X_i$.  Therefore, $\rmHN(X_i) =
\rmHN(X)|_{[0,\rk(X_i)]_\bfQ}$.  In particular, the breakpoints of the
HN polygon of $X$ are precisely the HN endpoints of the objects $X_i$.

\vskip 6pt

\noindent \mbox{\bf Proposition (7.3).} {\it Let $X$ be a nonzero
  object of $\calC$, with HN filtration $0 = X_0 \subsetneq X_1
  \subsetneq \cdots \subsetneq X_N = X$. Then:
  \begin{itemize}
  \item[(a)] For any a nonzero subobject $X'$ of $X$, $\deg(X') \leq
    \rmHN(X)(\rk(X'))$.  In other words, the HN endpoint of any
    subobject of $X$ lies on or below the HN polygon of $X$.
  \item[(b)] Moreover, if equality holds in (a), i.e.\ the HN endpoint
    of $X'$ lies on the HN polygon of $X$, then $X_{i-1} \subseteq X'
    \subseteq X_i$ for $i$ satisfying $\rk(X_{i-1}) \leq \rk(X') \leq
    \rk(X_i)$.
  \item[(c)] The HN polygon of $X$ is the upper convex hull of the
    points $(\rk(X'),\deg(X'))$, for $X'$ ranging over its nonzero
    subobjects.
  \item[(d)] For any nonzero subobject $X'$ of $X$, the HN polygon of
    $X'$ lies on or below the HN polygon of $X$.
  \end{itemize}}

\vskip 6pt

\noindent {\it Proof.} (a,b): At any point in time, it suffices to
replace $X'$ by $\wt{X}'$, so we may assume $X'$ is a strict subobject
of $X$.  Let us induct on the length $N$ of the HN filtration of $X$.
In the case $N=1$, $X$ is semistable.  But then $\mu(X') \leq \mu(X)$
and $\rk(X') \leq \rk(X)$ give (a); (b) is trivial.

Now assume $N>1$.  Note that if either $X' \subseteq X_1$ or $X_1
\subseteq X'$ holds, then we know (a,b) for $X'$.  In the first case,
one applies the base case above to $X_1$, and, in the second case, the
additivity of $\rk$ and $\deg$ allow us to reduce to the case of
$X'/X_1$ inside $X/X_1$, to which the inductive hypothesis applies.
By (6.2.a), we have $\mu(X') \leq \mu(X_1)$, and if equality holds
then $X' \subseteq X_1$; thus we assume also that $\mu(X') <
\mu(X_1)$.  Also, since $\rk(X' \cap X_1) = \rk(X_1)$ implies $X_1
\subseteq X'$, we can assume that $\rk(X' \cap X_1) < \rk(X_1)$.  And,
since $X'+X_1$ contains $X_1$, we know (a,b) for it.

Using, in order, (2.8), then the semistability of $X_1$, then that
$\mu(X') < \mu(X_1)$ but $\rk(X_1) > \rk(X' \cap X_1)$, we have
\begin{align*}
\rk(X'+X_1) \mu(X'+X_1)
& = \deg(X'+X_1) = \deg(X') + \deg(X_1) - \deg(X' \cap X_1) \\
& = \rk(X')\mu(X') + \rk(X_1)\mu(X_1)
  - \rk(X' \cap X_1)\mu(X' \cap X_1) \\
& \geq \rk(X')\mu(X') + [\rk(X_1)-\rk(X' \cap X_1)]\mu(X_1) \\
& > [\rk(X')+\rk(X_1)-\rk(X' \cap X_1)]\mu(X') \\
& = \rk(X'+X_1)\mu(X').
\end{align*}
Since $\rk(X'+X_1) > 0$, this implies $\mu(X') < \mu(X'+X_1)$.  This
inequality implies that the ray from the origin to the HN endpoint of
$X'$ lies strictly below the ray from the origin to the HN endpoint of
$X'+X_1$, and the endpoint of $X'$ clearly lies to the left of the
endpoint of $X'+X_1$.  By our knowledge of (a) for $X'+X_1$, we
condlude that the HN endpoint of $X'$ lies strictly beneath the HN
polygon of $X$, whence (a) and (b).

(c): Given (a), it suffices to see that the breakpoints of the HN
polygon are achieved.  But the breakpoints are precisely the HN
endpoints of the subobjects occurring in the HN filtration of $X$.

(d): Any subobject of $X'$ is a subobject of $X$, so the claim follows
from comparing the recipe of (c) as applied to $X$ and $X'$.
$\square$

\vskip 6pt

{\bf (7.4)} We give a common variant of the above.  A {\it semistable
filtration} of a nonzero object $X$ of $\calC$ is a filtration
$F_\bullet \cn 0 = F_0 \subsetneq F_1 \subsetneq \cdots \subsetneq F_M
= X$ by strict subobject such that each $F_i/F_{i-1}$ is semistable.
Its {\it polygon} $\rmP(F_\bullet)$ is obtained just like the HN
polygon is obtained from the HN filtration of $X$, except that one
must first sort the slopes of the $F_i/F_{i-1}$ into nonincreasing
order.

\vskip 6pt

\noindent \mbox{\bf Proposition (7.5).} {\it Let $X$ be a nonzero
  object of $\calC$, and let $0 = F_0 \subsetneq F_1 \subsetneq \cdots
  \subsetneq F_M$ be a semistable filtration.  Then $\rmHN(X) \leq
  \rmP(F_\bullet)$, and $\rmHN(X)(\rk(X)) = \rmP(F_\bullet)(\rk(X))$.}

\vskip 6pt

\noindent {\it Proof.} The coincidence of endpoints of the two
polygons follows immediately from the additivity of $\rk$ and $\deg$
over short exact sequences.  Let $0 = X_0 \subsetneq X_1 \subsetneq
\cdots \subsetneq X_N = X$ denote the HN filtration of $X$.  We first
show that, for each index $i = 1, \ldots, N$, one has
$\rmHN(X)(\rk(X_i)) \leq \rmP(F_\bullet)(\rk(X_i))$.

Note that, if we count the slopes of $F$ with multiplicities according
to the ranks, then for an integer $n$ with $0 \leq n \leq \rk(X)$,
$\rmP(F_\bullet)(n)$ is characterized as the sum of the $n$ largest
slopes.  Therefore, if we can find $n$ slopes of $F$ whose sum is
greater or equal to $C$, then we know that $\rmP(F_\bullet)(n) \geq
C$.  We apply this reasoning to $n = \rk(X_i)$ and $C = \deg(X_i)$.

Fix $i$.  For each $j = 1, \ldots, M$, the subobject $F_j \cap X_i
\subseteq X$ is strict by the last claim of (2.8).  Using the
strictness at $j-1$, we can form the quotient $(F_j \cap X_i)/(F_{j-1}
\cap X_i)$.  It is easy to check that the inclusion $F_j \cap X_i
\subseteq F_j$ induces a morphism $(F_j \cap X_i)/(F_{j-1} \cap X_i)
\to F_j/F_{j-1}$ that is an $F$-monomorphism, hence a monomorphism by
(2.3).  The semistability of $F_j/F_{j-1}$ gives the first inequality
in
\[
\mu\left(\frac{F_j \cap X_i}{F_{j-1} \cap X_i}\right)
  \leq \mu(F_j/F_{j-1})
\qquad \text{and} \qquad
d_j := \rk\left(\frac{F_j \cap X_i}{F_{j-1} \cap X_i}\right)
  \leq \rk(F_j/F_{j-1}).
\]
We claim that choosing the $j$th slope of $F$ with multiplicity $d_j$
works.  Indeed, the second inequality above shows that $d_j$ does not
exceed the multiplicity of the $j$th slope, and
\[
\sum_{j=1}^M d_j
= \sum_{j=1}^M [\rk(F_j \cap X_i) - \rk(F_{j-1} \cap X_i)] = \rk(X_i).
\]
Moreover,
\[
\sum_{j=1}^M d_j\mu(F_j/F_{j-1})
\geq \sum_{j=1}^M
    d_j\mu\left(\frac{F_j \cap X_i}{F_{j-1} \cap X_i}\right)
= \sum_{j=1}^M
    \left[\deg(F_j \cap X_i) - \deg(F_{j-1} \cap X_i)\right]
= \deg(X_i).
\]

Now, knowing that each $\rmHN(X)(\rk(X_i)) \leq
\rmP(F_\bullet)(\rk(X_i))$, we conclude.  For each $i$, consider the
restrictions of $\rmHN(X)$ and $\rmP(F_\bullet)$ to interval
$[\rk(X_{i-1}),\rk(X_i)]_\bfQ$.  The restriction of $\rmHN(X)$ is
linear, and the restriction of $\rmP(F_\bullet)$ is piecewise-linear,
concave down, and with endpoints above those of $\rmHN(X)$.
Therefore, we must have $\rmHN(X) \leq \rmP(F_\bullet)$ over this
interval.  Since these intervals cover $[0,\rk(X)]_\bfQ$, we are done.
$\square$

\vskip 6pt

\noindent \mbox{\bf Proposition (7.6).} {\it Let $0 \to X' \to X \to X'' \to
  0$ be a short strict exact sequence in $\calC$, with $X'$ and $X''$
  nonzero.  Then:
  \begin{itemize}
  \item[(a)] $\rmHN(X) \leq \rmHN(X' \oplus X'')$.
  \item[(b)] If all the slopes of $X'$ are greater than or equal to
    all the slopes of $X''$, then
    \[
    \rmHN(X' \oplus X'')(x) =
    \begin{cases}
      \rmHN(X')(x)
        & \text{if } 0 \leq x \leq \rk(X'), \text{ and} \\
      \rmHN(X'')(x-\rk(X')) + \deg(X')
        & \text{if } \rk(X') \leq x \leq \rk(X).
    \end{cases}
    \]
    In other words, the HN polygon of $X' \oplus X''$ is obtained by
    joining the HN polygons of $X'$ and $X''$.
  \end{itemize}}

\vskip 6pt

\noindent {\it Proof.} (a): By considering (7.3.c) applied to each of
$X$ and $X' \oplus X''$, it suffices to show that for any nonzero
strict subobject $Y$ of $X$, there exist subobjects $Y'$ of $X'$ and
$Y''$ of $X''$ with $\rk(Y) = \rk(Y')+\rk(Y'')$ and $\deg(Y) =
\deg(Y')+\deg(Y'')$, because then the subobject $Y' \oplus Y''$ would
bound the HN polygon of $X' \oplus X''$ above the HN polygon of $X$ at
$\rk(Y)$.  We may take $Y' = Y \cap X'$ and $Y'' = (Y+X')/X'$; the
computation of ranks is immediate, and the computation of degrees
follows from (2.10).

(b): The obvious filtration on $X' \oplus X''$ built from the HN
filtrations on $X'$ and $X''$ satisfies the characterizing properties
of the HN filtration.  (When $X'$ and $X''$ share a common slope, one
uses (4.8).)  The claim follows from this.  $\square$

\vskip 6pt

{\bf (7.7)} Given a sequence $0 \to X' \to X \to X'' \to 0$ as in
(7.6), it is natural to ask if conditions on the slopes of $X'$ and
$X''$ generally guarantee that the extension $X$ is split
(resp.\ nonsplit).  As far as I know, particular examples show that no
such general conditions exist: for any slope combination, one can find
both split and nonsplit extensions.  In other words, although
Harder--Narsimhan theory gives us information on Homs, such as in
(6.5), it does not give us information on higher $\Ext$s.

\vskip 12pt

\noindent {\bf \S8.  Examples: $\vphi$-modules.}

\vskip 6pt

We now give some variants of examples.

\vskip 6pt

\noindent \mbox{\bf Common Inputs (8.1).} All our examples will
involve:
\begin{itemize}
\item a B\'ezout domain $R$ with fraction field $K$, and
\item an injective ring homomorphism $\vphi \cn R \to R$.
\end{itemize}

Note that $\vphi$ extends uniquely to a ring homomorphism $K \to K$,
via the rule $\vphi(r/s) = \vphi(r)/\vphi(s)$.

\vskip 6pt

{\bf (8.2)} Because $R$ is a B\'ezout ring, any finitely generated,
torsion-free $R$-module is free.  In particular, a finitely generated
submodule $N$ of a finite free $R$-module $M$ is free.  Such $N$ is
moreover a strict subobject of $M$ in the exact category of finite
free $R$-modules if and only if it is a direct summand, if and only if
it is {\it saturated}:
\[
N = (K \otimes_R N) \cap M \text{ within } K \otimes_R M.
\]
In other words, for finite free $M$, the natural maps
\[
(8.2.1) \qquad
\left\{\begin{array}{@{}c@{}}
R\text{-direct summands} \\
\text{of }M
\end{array}\right\}
\subseteq
\left\{\begin{array}{@{}c@{}}
\text{saturated} \\
\text{$R$-submodules} \\
\text{of }M
\end{array}\right\}
\xrightarrow{K \otimes_R -}
\left\{\begin{array}{@{}c@{}}
\text{$K$-subspaces} \\
\text{of }K \otimes_R M
\end{array}\right\}
\]
are all bijections.

\vskip 6pt

{\bf (8.3)} If $M$ and $N$ are finite free $R$-modules and $f \cn M
\to N$ is an $R$-module map, we say that $f$ is an {\it isogeny} from
$M$ to $N$, written $f \cn M \stackrel\sim\dashrightarrow N$, if it
induces an isomorphism $f \cn K \otimes_R M \stackrel\sim\to K
\otimes_R N$.  This is the case if and only if $f$ is injective and
there exists nonzero $r \in R$ such that $rN \subseteq f(M)$.

\vskip 6pt

{\bf (8.4)} For an $R$-module $M$, define $\vphi^*M = R
\mathop{{}_\vphi\otimes_R} M$.  The subscripts on the $\otimes$-symbol
mean that $1 \otimes rm = \vphi(r) \otimes m$.  We consider this
object as an $R$-module using the left $\otimes$-factor.  To give an
$R$-linear map $\vphi^*M \to N$ is the same as to give a
$\vphi$-linear map $f' \cn M \to N$, that is, satisfying $f'(rm) =
\vphi(r)f'(m)$.  (Namely, $f'(m) = f(1 \otimes m)$ and $f(r \otimes m)
= rf'(m)$.)  If $M$ is finite free, then so is $\vphi^*M$, of the same
rank.  Analogous definitions and claims apply to $K$-vector spaces
instead of $R$-modules, compatibly with base change from $R$ to $K$.

\vskip 6pt

{\bf (8.4)} We define a {\it weak $\vphi$-module over $R$} to be a
finite free $R$-module $M$ equipped with an isogeny $\vphi_M \cn
\vphi^*M \stackrel\sim\dashrightarrow M$, and we let $\calC_\rmw$
denote the exact category of these (with $R$-linear maps respecting
the isogenies).  Such $(M,\vphi_M)$ is called a {\it $\vphi$-module
  over $R$} if $\vphi_M$ is an isomorphism, and the strictly full,
exact subcategory of these is written $\calC$.  Similarly, a {\it
  $\vphi$-module over $K$} is a finite-dimensional $K$-vector space
$M'$ equipped with an isomorphism $\vphi_{M'} \cn \vphi^*M'
\stackrel\sim\to M'$, and the abelian category of these will be
denoted $\calD$.

\vskip 6pt

{\bf (8.5)} If $M$ is an $R$-module and $k \geq 0$ an integer, we
write $M^{\wedge k}$ for the $k$th exterior power of $M$.  A map $f
\cn M \to N$ of $R$-modules gives rise to a map $f^{\wedge k} \cn
M^{\wedge k} \to N^{\wedge k}$ by the formula $f^{\wedge k}(m_1 \wedge
\cdots \wedge m_k) = f(m_1) \wedge \cdots \wedge f(m_k)$.  If $M$ is a
finite free $R$ module of rank $d$, then $M^{\wedge k}$ is a finite
free $R$-module of rank $\binom d k$.  If $f$ is an isogeny between
free $R$-modules of rank $d$, then $f^{\wedge d}$ is an isogeny
between free $R$-modules of rank $1$.

In particular, if $(M,\vphi_M)$ is a weak $\vphi$-module over $R$ of
rank $d$, then
\[
\vphi_{M^{\wedge d}} \cn
\vphi^*(M^{\wedge d})
\cong (\vphi^*M)^{\wedge d}
\xrightarrow{\vphi_M^{\wedge d}} M^{\wedge d}
\]
defines a weak $\vphi$-module $(M^{\wedge d},\vphi_{M^{\wedge d}})$
over $R$ of rank $1$.  Choosing a basis $v \in M^{\wedge d}$, we let
$\det(\vphi_M) \in R$ be the nonzero element satisfying
$\vphi_{M^{\wedge d}}(1 \otimes v) = \det(\vphi_M) v$.  It is defined
up to multiplication by $\vphi(r)/r$, for $r \in R^\times$.  Note that
$(M,\vphi_M)$ is a $\vphi$-module if and only if $(M^{\wedge
  d},\vphi_{M^{\wedge d}})$ is, if and only if $\det(\vphi_M) \in
R^\times$.

Analogous definitions and claims apply to $K$-vector spaces
instead of $R$-modules, compatibly with base change from $R$ to $K$.

\vskip 6pt

{\bf (8.6)} We apply HN theory to:
\begin{gather*}
\calC = \calC_\rmw \text{ or } \calC, \qquad
\calD = \calD, \qquad
F = K \otimes_R -, \qquad
\rk = \dim_K,
\end{gather*}
and $V$ and $\deg$ to be specified below.  Let us check the axioms
that concern these data.

\vskip 6pt

{\bf (8.7)} We check the axiom for $F$.  Faithfulness and exactness
are clear, so we treat the remaining claim.  Fix an object
$(M,\vphi_M)$ of $\calC_\rmw$ or $\calC$, and for brevity let
$F(M,\vphi_M) = (M',\vphi_{M'})$.  A subobject of $(M,\vphi_M)$ in
$\calC_\rmw$ is uniquely determined by an $R$-submodule $N$ of $M$
such that $\vphi_M(\vphi^*N) \subseteq N$, and this subobject is
strict if and only if $N$ is saturated in $M$.  Similarly, a subobject
of $(M',\vphi_{M'})$ in $\calD$ is uniquely determined by a
$K$-subspace $N'$ of $M'$ such that $\vphi_{M'}(\vphi^*N') \subseteq
N'$ (whence comparing dimensions gives $\vphi_{M'}(\vphi^*N') = N'$).
For a saturated $R$-submodule $N$ of $M$ that corresponds to the
$K$-subspace $N'$ of $M'$ under (8.3.1), so that $N' = K \otimes_R N$
and $N = N' \cap M$, we must show that $\vphi_M(\vphi^*N) \subseteq N$
if and only if $\vphi_{M'}(\vphi^*N') \subseteq N'$.  Indeed, if the
latter holds then
\[
\vphi_M(\vphi^*N)
\subseteq \vphi_{M'}(\vphi^*N') \cap \vphi_M(\vphi^*M)
\subseteq N' \cap M = N
\]
as was desired, and the reverse implication is clear.

When $(M,\vphi_M)$ belongs to $\calC$, we must moreover show that when
$N'$ is a $\vphi_{M'}$-stable $K$-subspace of $M'$, the saturated
$R$-submodule $N = N' \cap M$ satisfies $\vphi_M(\vphi^*N) = N$.  In
fact, the above argument shows that $N$ inherits a weak $\vphi$-module
structure, and we also get a weak $\vphi$-module structure on $M/N$.
Then one has $\det(\vphi_M) = \det(\vphi_N)\det(\vphi_{M/N})$, and
$\det(\vphi_M)$ is a unit, so $\det(\vphi_N)$ is also a unit.

\vskip 6pt

{\bf (8.8)} The axiom for $\rk$ is trivial.

\vskip 6pt

\noindent \mbox{\bf First Variant (8.9).} We let $V$ be a totally
ordered abelian group (written additively), and $\ord \cn K^\times \to
V$ a homomorphism, satisfying:
\begin{itemize}
\item for all nonzero $r \in R$ one has $\ord(\vphi(r)/r) \geq 0$,
  with equality if and only if $r \in R^\times$.
\end{itemize}
Then we take $\calC = \calC_\rmw$, and $\deg(M,\vphi_M) =
-\ord(\det(\vphi_M))$.

\vskip 6pt

{\bf (8.10)} We verify the axiom for $\calC_\rmw$ and this $\deg$.
The additivity of $\deg$ follows from multiplicativity of the
determinant over short exact sequences, and that $\ord$ is a
homomorphism.  To verify the other claim, by passing to top exterior
powers we may reduce to the case of an $F$-isomorphism $f \cn
(M,\vphi_M) \to (N,\vphi_N)$ with $M$ and $N$ of rank $1$; we identify
$M$ with its image in $N$ under $f$.  Choose a basis $v$ of $N$ and a
nonzero $r \in R$ such that $M$ is generated by $rv$.  We have
$\vphi_N(1 \otimes v) = \det(\vphi_N) v$ and
\[
\det(\vphi_M) \cdot rv
= \vphi_M(1 \otimes rv)
= \vphi(r)\vphi_N(1 \otimes v)
= \vphi(r) \det(\vphi_N) v
= (\vphi(r)/r)\det(\vphi_N) \cdot rv,
\]
so $\deg(M) = \deg(N) - \ord(\vphi(r)/r)$.  The claim is now seen to
be equivalent to our axiom for $\ord$.

\vskip 6pt

\noindent \mbox{\bf Example (8.11).}  Suppose $R$ is a valuation ring
(and therefore a B\'ezout domain).  The group $V = K^\times/R^\times$,
called the value group of $R$, is totally ordered by $a \leq b$ if and
only if $aR \supseteq bR$ within $K$.  (We ignore that the operation
of $V$ is by default written multiplicatively.)  The projection map
$K^\times \to V$, called the valuation, will serve as our $\ord$.  For
$r \in K^\times$ one has $\ord(r) \geq 0$ if and only if $r \in R$,
with equality if and only if $r \in R^\times$.  For these data, our
axiom becomes the following requirement on an injective ring
homomorphism $\vphi \cn R \to R$: for all nonzero $r \in R$, one has
$\vphi(r) \in rR$, and $\vphi(r)R = rR$ implies $r \in R^\times$.

Let us explicitly construct such an $R$ and $\vphi$.  Let $k$ be a
field with an automorphism $\sigma \cn k \stackrel\sim\to k$, $X$ an
indeterminate, and $q \geq 2$ an intger.  Let $R = k\bbr X$, and
$\vphi \cn \sum_{n \geq 0} a_n X^n \mapsto \sum_{n \geq 0} \sigma(a_n)
X^{qn}$.  Then $R$ is a vlauation ring for the $X$-adic valuation and
the identity $\ord \circ \vphi = q \ord$ implies the required axiom.

\vskip 6pt

\noindent \mbox{\bf Second Variant (8.12).}  Let $V$ be a totally
ordered abelian group (written additively), and $\ord \cn R^\times \to
V$ a homomorphism, satisfying:
\begin{itemize}
\item for all nonzero $r \in R$ such that $\vphi(r)/r \in R^\times$,
  one has $\ord(\vphi(r)/r) \geq 0$, with equality if and only if $r
  \in R^\times$.
\end{itemize}
Then we take $\calC = \calC$ and $\deg(M,\vphi_M) =
-\ord(\det(\vphi_M))$.

\vskip 6pt

{\bf (8.13)} The axiom for $\calC$ and this $\deg$ are verified just
as in (8.10).  (In the notation there, the fact that both
$\det(\vphi_M)$ and $\det(\vphi_N)$ are units implies $\vphi(r)/r$ is
also a unit, so $\deg(\vphi(r)/r)$ is defined.)

\vskip 6pt

\noindent \mbox{\bf Example (8.14).}  Fix a prime $p$, and let $R
\subset \bfQp\bbr{X,X\inv}$ be the Robba ring.  This consists of
Laurent series $f = \sum_{n \in \bfZ} a_nX^n$ convergent in some
$p$-adic annulus of the form $\ep \leq |X| < 1$, where $\ep$ is
allowed to depend on $f$.  There is a unique continuous ring
endomorphism $\vphi$ that is the identity on coefficients and sends $X
\mapsto (1+X)^p-1$.  (The $\ep$ of convergence for $\vphi(f)$ may be
closer to $1$ than the $\ep$ of convergence for $f$.)

Some nontrivial facts include that $R$ is a B\'ezout ring, and the
following description of its units.  Consider the subring $E \subset
R$ of series $f$ that, after perhaps adjusting $\ep$, are not only
convergent but bounded on the annulus $\ep \leq |X| < 1$.  Then the
inclusion $E^\times \subseteq R^\times$ is a bijection.  Moreover, $E$
is a (Henselian) discretely valued field with uniformizer $p$, and its
valuation $\ord \cn R^\times = E^\times \to \bfZ$ satisfies the
hypothesis of (8.12).  This example is the subject of [K].

\vskip 6pt

\noindent \mbox{\bf Non-Example (8.15).} Let $k$ be a perfect field of
characteristic $p>0$.  We take $V=\bfZ$, $R = W(k)$ the Witt vectors
of $k$, $\vphi$ to be the automorphism induced from the $p$th power
Frobenius map on $k$, and $\ord$ to be the $p$-adic valuation.  In
this setting, $\calC_\rmw$ is the category of so-called {\it
  $F$-crystals} over $k$.  Then $\ord \circ \vphi = \ord$, so neither
variant of this section applies.  In fact, with our choice of $\calD$
and $F$ the axiom for $\deg$ always fails: for any nonzero object
$(M,\vphi_M)$, the inclusion $pM \subseteq M$ is an $F$-isomorphism
but not an isomorphism, whereas the two isomorphic objects $pM \cong
M$ would have the same degree.

Nonetheless, there is an adaptation of HN theory to $\calC_\rmw$
resembling the variant (8.9), using the $p$-adic valuation for
$\ord_p$ to define $\deg(M,\vphi_M) = -\ord_p(\det(\vphi_M))$, but
replacing the use of $F \cn \calC_\rmw \to \calD$ with other means.
This, and more, is the subject of the Dieudonn\'e--Manin Theorem.

\vskip 18pt

\noindent {\bf References.}

\begin{itemize}
\item[{[F]}] L.~Fargues, La filtration de Harder--Narasimhan des
  sch\'emas en groupes finis et plats.  {\it J. Reine Angew.\ Math.}
  {\bf 645} (2010), 1--39.
\item[{[K]}] K.~S.~Kedlaya, Slope filtrations revisited.  {\it
  Doc.\ Math.} {\bf 10} (2005), 447--525.
\end{itemize}

\end{document}